\documentclass[12pt,a4paper,reqno]{amsart}
\usepackage{amsmath}
\usepackage{amsfonts}
\usepackage{amssymb}
\numberwithin{equation}{section}

\theoremstyle{plain}

\newtheorem{theorem}[subsection]{Theorem}
\newtheorem{proposition}[subsection]{Proposition}
\newtheorem{lemma}[subsection]{Lemma}
\newtheorem{corollary}[subsection]{Corollary}

\newtheorem*{weyl-theorem}{Weyl's Inequality}

\theoremstyle{definition}

\newtheorem{definition}[subsection]{Definition}

\theoremstyle{remark}

\newtheorem*{example}{Example}

\newtheorem*{remark}{Remark}

\renewcommand{\leq}{\leqslant}
\renewcommand{\geq}{\geqslant}

\newsavebox{\proofbox}
\savebox{\proofbox}{\begin{picture}(7,7)%
  \put(0,0){\framebox(7,7){}}\end{picture}}
\def\boxeq{\tag*{\usebox{\proofbox}}}

\newcommand{\md}[1]{\ensuremath{(\mbox{mod}\, #1)}}
\def\proof{\noindent\textit{Proof. }}
\def\endproof{\hfill{\usebox{\proofbox}}}
\def\E{\mathbb{E}}

\def\Z{\mathbb{Z}}
\def\R{\mathbb{R}}

\def\C{\mathbb{C}}
\def\P{\mathbb{P}}

\def\B{\mathcal{B}}

\def\triv{{\operatorname{triv}}}

\def\eps{\varepsilon}

\begin{document}

\title[A new bound for $r_4(N)$]{New bounds for Szemer\'edi's theorem, II: A new bound for $r_4(N)$}

\author{Ben Green}
\address{Centre for Mathematical Sciences, Wilberforce Road, Cambridge CB3 0WA, England.}
\email{b.j.green@dpmms.cam.ac.uk}

\author{Terence Tao}
\address{Department of Mathematics, UCLA, Los Angeles CA 90095-1555, USA.
}
\email{tao@math.ucla.edu}

\thanks{The first author is a Clay Research Fellow, and is pleased to acknowledge the support of the Clay Mathematics Institute. Some of this work was carried out while he was on a long-term visit to MIT. The second author is supported by a grant from the Packard Foundation.}

\begin{abstract} 
Define $r_4(N)$ to be the largest cardinality of a set $A \subseteq \{1,\dots,N\}$ which does not contain four elements in arithmetic progression. In 1998 Gowers proved that
\[ r_4(N) \ll N(\log \log N)^{-c}\] for some absolute constant $c> 0$. In this paper (part II of a series) we improve this to
\[ r_4(N) \ll N e^{-c\sqrt{\log\log N}}.\]
In part III of the series we will use a more elaborate argument to improve this to 
\[ r_4(N) \ll N(\log N)^{-c}.\]
\end{abstract}

\maketitle
\begin{center}
\emph{To Klaus Roth on his 80th birthday}
\end{center}

\section{Introduction}

\textsc{notational convention.} Throughout the paper the letters $c,C$ will denote absolute constants which could be specified explicitly if desired. These constants will generally satisfy $0 < c \ll 1 \ll C$. Different instances of the notation, even on the same line, will typically denote different constants. Occasionally we will want to fix a constant for the duration of an argument; such constants will be subscripted as $C_0,C_1$ and so on. Any implied constants in the $O$- or $\ll$ notations will depend only on any subscripted variables. Thus if we say that $f(N) = O_{\delta}(N)$ we mean that there is a constant $F(\delta)$ such that $f(N) \leq F(\alpha) N$ for all $N$. The absence of any subscripted variables should be taken to mean that the implied constant is absolute.

Let $N$ be a large positive integer, and let $k \geq 3$ be fixed. We define $r_k(N)$ to be the largest cardinality of a set $A \subseteq [N] = \{1,\dots,N\}$ which does not contain $k$ distinct elements in arithmetic progression. 

Klaus Roth proved in 1953 \cite{roth} that 
\[ r_3(N) \ll N(\log \log N)^{-1}.\]
In particular, $r_3(N) = o(N)$.
Since Szemer\'edi's 1969 proof that $r_4(N) = o(N)$ \cite{szemeredi-4}, and his later proof \cite{szemeredi} that $r_k(N) = o_k(N)$ for $k \geq 5$, it has been natural to ask for similarly effective bounds for these quantities. A first attempt in this direction was made by Roth in \cite{roth-4}, who provided a new proof that $r_4(N) = o(N)$. A major breakthrough was made by in 1998 by Gowers \cite{gowers-4-aps,gowers-long-aps}, who obtained the bound
\[ r_k(N) \ll N(\log \log N)^{-c_k}\] for each $k \geq 4$. 

In the meantime, there has been progress on $r_3(N)$. Szemer\'edi (unpublished) obtained the bound 
\begin{equation}\label{eq1.1} r_3(N) \ll Ne^{-c\sqrt{\log \log N}},\end{equation} and shortly thereafter
Heath-Brown \cite{heath} and Szemer\'edi \cite{szem-3ap} independently obtained the bound
\begin{equation}\label{eq1.2} r_3(N) \ll N(\log N)^{-c}.\end{equation}
More recently Bourgain \cite{Bou} found the best bound currently known, namely
\[ r_3(N) \ll N(\log \log N/\log N)^{1/2}.\]
Part I of this series of papers \cite{green-tao-szem} may be consulted for a more extensive discussion of the history of the problem. Our objective in this series is to bring our knowledge of $r_4$ more closely into line with the best known bounds for $r_3$. In \cite{green-tao-szem} this was achieved in the so-called \emph{finite field model}, in which $[N]$ is replaced by a vector space $\mathbb{F}_p^n$ over a finite field. In this paper we instead study subsets of $[N]$ itself, and obtain the analogue of Szemer\'edi's unpublished bound \eqref{eq1.1} for $r_4(N)$. 

\begin{theorem}[Main theorem]\label{main}  For all large integers $N$ we have
\[ r_4(N) \ll N e^{- c \sqrt{\log \log N}}.\]
\end{theorem}

In part III of the series we will obtain the analogue of the superior bound \eqref{eq1.2}. The argument will, however, be substantially more technical.

Let us conclude this introduction by mentioning that the best known \emph{lower} bound for $r_4(N)$ is essentially the same as that for $r_3(N)$, namely Behrend's 1946 bound \cite{behrend}
\[ r_4(N) \geq r_3(N) \gg Ne^{-c\sqrt{\log N}}.\]
Somewhat better bounds of shape
\[ r_k(N) \gg Ne^{-(\log N)^{c_k}}\]
are known for much larger $k$: see  \cite{laba,rankin} for details.

We now briefly outline the proof of Theorem \ref{main}.  As with all previous papers obtaining quantitative bounds for $r_k(N)$, we use the \textit{density increment strategy} of Roth, a detailed discussion of which may be found in \cite{green-fin-field}.  
The key is to obtain a dichotomy of the following form.

\begin{proposition}[Lack of progressions implies density increment]\label{densinc-general}  Let $N$ be a large integer, let $\delta \in (0,1)$, and suppose that $A \subseteq [N]$ has $|A| \geq \delta N$ and contains no progressions of length $4$.  Assume that we have the largeness condition $N \geq F(\delta)$ for some explicit function $F$.
Then there exists an arithmetic progression $P \subseteq [N]$ of length at least $f(N,\delta)$ on which we have the density increment
\[ \frac{|A \cap P|}{|P|} \geq \delta + \sigma(\delta).\]
Here $f(N,\delta) > 0$ is an explicit function which goes to $\infty$ as $N \to \infty$ for each fixed $\delta$, and $\sigma(\delta) > 0$ is an explicit positive quantity depending only on $\delta$.
\end{proposition}

Any proposition of this type will imply, by iteration, a nontrivial upper bound on $r_4(N)$, with the precise bound depending on the functions $F(\,)$, $f(\, , \,)$, and $\sigma(\,)$. For an actual calculation of a bound (on $r_4(\mathbb{F}_5^n)$) using this strategy, part I of this series may be consulted.

If one desires a good bound it is of particular interest to get $f(N,\delta)$ and $\sigma(\delta)$ as large as possible. The function $F(\delta)$ plays a much less significant r\^ole and, at least for the purposes of a motivating discussion, may be ignored.  Gowers' proof that $r_4(N) \ll N(\log \log N)^{-c}$ proceeds by establishing Proposition \ref{densinc-general} with
$f(N,\delta) \gg N^{c\delta^C}$ and $\sigma(\delta) \gg \delta^C$.
The main advance in our paper is to improve the density increment bound to $\sigma(\delta) \gg \delta$. This has the effect of reducing the number of iterations of Proposition \ref{densinc-general} that are required from $C\delta^{-C}$ to $C \log(1/\delta)$. Here is a more precise statement of what we shall prove.

\begin{proposition}[Lack of progressions implies density increment]\label{densinc-specific}  Let $\delta > 0$, and suppose that $N \geq e^{C \delta^{-C}}$. Let $A$ be a subset of $[N]$ with $|A| \geq \delta N$ such that
$A$ contains no progressions of length $4$.  Then there exists an arithmetic 
progression $P$ in $[N]$ of length $|P| \gg N^{c \delta^C}$ such that we have the density increment
\[ \frac{|A \cap P|}{|P|} \geq (1+c) \delta.\]
\end{proposition}

Let us now quickly show how this implies Theorem \ref{main}.

\emph{Deduction of Theorem \ref{main} from Proposition \ref{densinc-specific}.}
Suppose that $A \subseteq [N]$ has size $\delta N$, and that it does not contain a 4-term progression. We perform an iteration. At the $i$th step of this iteration we will have a set $A_i \subseteq \{1,\dots,N_i\}$ with size $\delta_i N$.  This set will be a linearly rescaled version of a subset of $A$, and so it too does not contain a progression of length 4. Set $A_0 := A$, $N_0 := N$ and $\delta_0 := \delta$. Now Proposition \ref{densinc-specific} tells us that either
\begin{equation}\label{bottom-out}
N_i \leq e^{C\delta_i^{-C}}
\end{equation}
or else the iteration proceeds and it is possible to choose $N_{i+1}, \delta_{i+1}$ and $A_{i+1}$ such that\[ N_{i+1} \gg N_i^{c\delta_i^C}\]
and \[ \delta_{i+1} \geq (1 + c)\delta_i.\]
Now as long as the iteration continues we must have $\delta_i \leq 1$, and so after $K \leq C \log(1/\delta)$ iterations the condition \eqref{bottom-out} must be satisfied. At this point we have
\[ N_K \gg N^{(c\delta^C)^{C\log(1/\delta)}},\] and so we derive the inequality
\[ N^{(c\delta^C)^{C\log(1/\delta)}} \leq e^{C \delta^{-C}}.\]
After a small amount of rearrangement this leads to the claimed bound
\begin{equation}\boxeq r_4(N) \ll N e^{-c\sqrt{\log \log N}}.\end{equation}

It thus remains to establish Proposition \ref{densinc-specific}.  Our starting point is our earlier paper \cite{gt-inverseu3}, which built upon the original paper of Gowers \cite{gowers-4-aps} to provide an \emph{inverse $U^3$ theorem} which, among other things, already implies Gowers' bound $r_4(N) \ll N(\log \log N)^{-c}$.  This inverse theorem will be stated properly in later sections, but roughly speaking if $A \subset [N]$ had size $|A| \geq \delta N$
and had no progressions of length $4$, then $A$ would have significant correlation with a certain ``local quadratic phase function''. An example of such a function is $n \mapsto e^{2\pi i \alpha n^2}$, though this is not the most general example; the reader may wish to consult the surveys \cite{green-icm,tao-coates} for further discussion.

This correlation implies that $A$ has a significant density increment (comparable to $\delta^C$) on a ``quadratic Bohr set'', that is to say an approximate level set of a local quadratic phase function. Such a set has size $\sim \delta^C N$. The next step in \cite{gowers-4-aps} is then that of \emph{linearisation}, in which the Bohr set is partitioned into arithmetic progressions. By the pigeonhole principle $A$ also has a density increment of $\sim \delta^C$ on one of these progressions. It turns out that the linearisation can be achieved with progressions of size $\gg N^{c\delta^C}$ which, as mentioned earlier, is sufficient to give Gowers' bound. We remark that a similar linearisation step already appears in the earlier work of Roth \cite{roth} (cf. \cite[Lemma 2.3]{gowers-long-aps}).

The main cost in this scheme lies in the linearisation step, which forces one to pass from an object of size $N$ to an object of size only $N^{c\delta^C}$.
To improve upon this scheme we borrow an idea of Heath-Brown and Szemer\'edi \cite{heath,szem-3ap} from the $k=3$ case. Instead of finding a quadratic phase function which correlates with $A$ and then linearizing, one adopts a more patient stance and first collects \emph{several} quadratic phase functions. In this way a more substantial density increment of $c\delta$ can be obtained. Only after this is done do we linearise. This procedure of linearizing several quadratic phase functions at once turns out not to be as costly as one might think, and in any case it need only be done $O(\log (1/\delta))$ times due to the size of the density increment.

In part III of the series we will show that it is possible to be more efficient still, by extracting additional gains either from the density increment or from the length of the progression on which the increment is obtained. This was carried out in the finite field setting in \cite{green-tao-szem}.

We have mentioned, albeit briefly, the so-called finite field model: the survey \cite{green-fin-field} may be consulted for more information. The advantage of working in $\mathbb{F}_p^n$ as opposed to the cyclic group $\Z/N\Z$ (which serves as a model for $[N]$) is the availability of \emph{subspaces}. In $\Z/N\Z$, and in other abelian groups $G$, one must make do with the notion of \emph{Bohr sets}, which may be thought of as approximate subspaces. There are various technical issues involved in dealing with these, as we shall see later on.

\begin{remark} It is quite likely that the methods here combine with those in \cite{green-tao-szem} extend to general finite abelian groups $G$; thus if $r_4(G)$ denotes the largest cardinality $|A|$ of a set $A \subset G$ without any arithmetic progressions of length $4$, a slight elaboration of the arguments here should establish $r_4(G) \ll |G| e^{- c \sqrt{\log \log |G|} }$ for all large $|G|$. We will however not pursue this matter here.
\end{remark}

\section{General notation}\label{notation-sec}

Let $A$ be a finite non-empty set and let $f: A \to \C$ be a function. It is convenient, so as to avoid having to contend with normalising factors, to use the expectation notation
\[\E_A(f) = \E_{x \in A} f(x) := \frac{1}{|A|} \sum_{x \in A} f(x).\]
More complex expressions such as $\E_{x \in A, y \in B} f(x,y)$ are similarly defined.  We also define the $L^p$ norms
$$ \|f\|_{L^p(A)} := (\E_A|f|^p)^{1/p}$$
for $1 \leq p < \infty$, with the usual convention $\|f\|_{L^\infty(A)} := \sup_{x \in A} |f(x)|$.  
We say that $f$ is \emph{$1$-bounded} if $\|f\|_{L^\infty(A)} \leq 1$.

If $A,B$ are finite sets with $B$ non-empty, we write $\P_B(A) := \frac{|A \cap B|}{|B|}$ for the density of $A$ in $B$. If $A$ lies in some ambient space $X$ (for example a group) we use $1_A: X \to \R$ to denote the 
indicator function of $A$, that is to say $1_A(x) = 1$ when $x \in A$ 
and $1_A(x) = 0$ otherwise. We also write $1_{x \in A}$ for $1_A(x)$.  Thus for instance $\P_B(A) = \E_B(1_A)$ for all 
non-empty $B \subseteq X$.

\section{The form $\Lambda$ and the $U^3(\Z/p\Z)$ norm}\label{lambda-sec}

We now begin the proof of Proposition \ref{densinc-specific}.  It will be convenient to work in a cyclic group $\Z/p\Z$ of
large prime order $p$ rather than on the interval $[N]$.  On this cyclic group $\Z/p\Z$, we introduce the quadrilinear form $\Lambda(f_0,f_1,f_2,f_3)$, defined for four functions
$f_j: \Z/p\Z \to \C$ by
\[ \Lambda(f_0,f_1,f_2,f_3) := \E_{x,h \in \Z/p\Z} f_0(x) f_1(x+h) f_2(x+2h) f_3(x+3h).\]
This form is clearly pertinent to the task of counting progressions of length $4$, and has appeared in many previous papers on this subject.
One can quickly deduce Proposition \ref{densinc-specific}, and hence Theorem \ref{main}, from the following claim.

\begin{theorem}[Anomalous number of AP4s implies density increment]\label{dens-inc}
Let $p$ be a large prime, let $N$ be an integer between $p/8$ and $p/4$, and let $f: \Z/p\Z \to \R$ be a $1$-bounded non-negative
function which vanishes outside of $[N]$.  Set $\delta := \E_{[N]}(f)$.  Suppose that 
\begin{equation}\label{plarge}
p \gg \exp( C \delta^{-C} )
\end{equation}
for some suitably large
absolute constant $C$, and suppose that
\begin{equation}\label{law}
|\Lambda(f,f,f,f) - \Lambda(\delta 1_{[N]}, \delta 1_{[N]}, \delta 1_{[N]}, \delta 1_{[N]})| \gg \delta^4.
\end{equation}
Then we can find an arithmetic progression $P$ in $[N]$ obeying the length bound 
\begin{equation}\label{qinc}
|P| \gg p^{c \delta^C}
\end{equation}
and the density increment bound
\begin{equation}\label{pinc}
 \E_P(f) \geq (1+c) \delta
\end{equation}
for some $c, C > 0$.
\end{theorem}

\begin{remark}
Strictly speaking, there could be two different notions of an arithmetic progression in $[N]$, one arising from its embedding into the integers $\Z$, and the other arising from its embedding into the cyclic group $\Z/p\Z$.  However, because $N < p/4$, it is easy to see that the two concepts are equivalent; the interval $[N]$ is too short to contain a progression that somehow ``wraps around'' $p$.  (To use some jargon, the two representations of $[N]$ are \emph{Freiman isomorphic} of order $2$, which is sufficient to preserve the concept of an arithmetic progression; see for instance \cite{taovu-book}.)
\end{remark}

\emph{Proof that Proposition \ref{densinc-specific} implies Theorem \ref{dens-inc}.}  By increasing $\delta$ if necessary we may
assume that $|A| = \delta N$.  Choose a prime $p$ between $4N$ and $8N$ (this is of course possible by Bertrand's Postulate)
and take $f := 1_A$, thought of as a function on $\Z/p\Z$.  Since $A$ has no progressions of length $4$ we easily see that
\[ \Lambda(f,f,f,f)= O( 1/p ),\]
whilst the fact that there are $\frac{1}{6}N^2 (1 + o(1))$ four-term progressions in $[N]$ implies that
\[ \Lambda(\delta 1_{[N]}, \delta 1_{[N]}, \delta 1_{[N]}, \delta 1_{[N]}) \gg \delta^4.\]
Since we are taking $p$ to be large, we conclude \eqref{law}.  Applying Theorem \ref{dens-inc}, we can find a progression $P$ in $\Z/p\Z$
obeying \eqref{qinc} and \eqref{pinc}, and this suffices for our needs.  \endproof

It remains to prove Theorem \ref{dens-inc}.  For the rest of the paper we fix $p$ to be a large prime.
To be able to exploit the hypothesis \eqref{law}, we will need to show that $\Lambda$ is controlled by either of two norms (when restricted
to $1$-bounded functions).  The first is the $L^1$ norm.

\begin{lemma}[$L^1$ controls $\Lambda$]\label{l1-lambda4} Let $f, g: \Z/p\Z \to \C$ be uniformly bounded by some $\alpha > 0$.
Then we have
\[ |\Lambda(f,f,f,f) - \Lambda(g,g,g,g)| \leq 4 \alpha^3 \|f-g\|_{L^1(\Z/p\Z)}.\]
\end{lemma}
\proof
Since $\Lambda$ is quadrilinear we have
\begin{align}\nonumber \Lambda(f,f,f,& f) - \Lambda(g,g,g,g) \\&= \Lambda(f-g,f,f,f) + \Lambda(g,f-g,f,f) + \Lambda(g,g,f-g,f) + \Lambda(g,g,g,f-g).\label{telescope}\end{align}
The result now follows on applying the triangle inequality and the easily checked bound
\begin{equation}\label{fairly-trivial} |\Lambda(f_1,f_2,f_3,f_4)| \leq \Vert f_j \Vert_1 \sup_{i = 1,\dots,4} \Vert f_i \Vert_{\infty}^3,\end{equation}
valid for $j = 1,\dots,4$.\endproof

The second norm that controls $\Lambda$ is the \emph{Gowers $U^3$-norm} $\|f\|_{U^3(\Z/p\Z)}$ of a function $f: \Z/p\Z \to \C$, defined as
\begin{eqnarray*} \Vert f \Vert_{U^3(\Z/p\Z)}^8 & := & \E_{x,h_1,h_2,h_3 \in \Z/p\Z}(f(x)\overline{f(x+h_1)f(x+h_2)f(x+h_3)}f(x+h_1+h_2) \times \\
& & \qquad\qquad\qquad\qquad \times f(x+h_2+h_3)f(x+h_1+h_3) \overline{f(x+h_1+h_2+h_3)}).\end{eqnarray*}
This norm was introduced in \cite{gowers-4-aps,gowers-long-aps} and
studied further in such papers as \cite{green-tao-primes,gt-inverseu3,green-tao-szem,taovu-book}.  As shown in \cite{gowers-long-aps} it is indeed a norm on $\Z/p\Z$, but we will not need to know this here.  In fact we only require two facts about the $U^3$-norm. One of these facts is an \emph{inverse theorem}, which will be the subject of the next section.
The other is that the $U^3$-norm controls $\Lambda$.

\begin{lemma}[$U^3$ controls $\Lambda$]\label{u3-lambda4} 
Let $f, g: \Z/p\Z \to \C$ be $1$-bounded functions on an affine space $\Z/p\Z$.  Then we have
$$ |\Lambda(f,f,f,f) - \Lambda(g,g,g,g)| \leq 4 \|f-g\|_{U^3(\Z/p\Z)}.$$
\end{lemma}
\proof We employ the same telescoping identity \eqref{telescope} that we use to prove Lemma \ref{l1-lambda4}. In place of the fairly trivial bound \eqref{fairly-trivial} we instead apply the \emph{Generalized von Neumann theorem}, which in this setting states that
\[ |\Lambda(f_1,f_2,f_3,f_4)| \leq \Vert f_j \Vert_{U^3(\Z/p\Z)}\]
for $j = 1,\dots,4$. This result is proved using three applications of the Cauchy-Schwarz inequality: the details are given very explicitly in \cite[Proposition 1.11]{green-montreal}.\endproof

In \cite{gowers-4-aps,gowers-long-aps} one applied Lemma \ref{u3-lambda4} directly to \eqref{law} in order to obtain the lower bound 
$\| f - \delta 1_{[N]} \|_{U^3(\Z/p\Z)} \gg \delta^3$.  This ultimately led to a density increment of $\gg \delta^C$ for $f$ on some progression. The resulting iteration scheme thus proceeds for $\gg \delta^{-C}$ steps, which is too long for our purposes. Our approach is to develop a so-called \emph{Koopman-von Neumann structure theorem}, which introduces an intermediate approximant $\E(f|\B)$ between $f$ and $\delta 1_{[N]}$. 

\section{The inverse $U^3(\Z/p\Z)$ theorem}\label{inverse-sec}

We now come to the second fact concerning the $U^3(\Z/p\Z)$-norm that we shall need, namely the \emph{inverse $U^3$-theorem}. This is one of the main results of \cite{gt-inverseu3}.  There are three (equivalent) formulations of this inverse theorem: one involving locally quadratic phase functions, one involving generalized quadratic phases, and one involving $2$-step nilsequences.  Our argument would work with the first two of these but not the third (cf. \cite[Theorem 12.7]{gt-inverseu3}), which has rather weaker bounds. We use the first formulation involving locally quadratic phases. This is in a sense the most basic form of the inverse theorem for $U^3(\Z/p\Z)$, since in \cite{gt-inverseu3} the other variants are all derived from it. To describe the result we need some notation.

\begin{definition}[Bohr sets]\label{bohr-set}  
Let $S \subseteq \Z/p\Z$, and let $\rho \in (0,1)$ be a parameter. We define the \emph{\textup{(}centred\textup{)} Bohr set} $B(S,\rho) \subseteq \Z/p\Z$ to be the set
\[ B(S,\rho) := \{ x \in \Z/p\Z : \|\xi x / p\|_{\R/\Z} < \rho \},\]
where $\|x\|_{\R/\Z}$ denotes the distance from $x$ to the nearest integer.  More generally, if $\alpha = (\alpha_\xi)_{\xi \in S}$ is any element in the $|S|$-dimensional torus $(\R/\Z)^S$, then we write $B_\alpha(S,\rho)$ for the \emph{uncentred Bohr set}
\[ B_\alpha(S,\rho) := \{ x \in \Z/p\Z : \|\xi x / p - \alpha_\xi \|_{\R/\Z} < \rho \}.\]
We refer to $|S|$ as the \emph{rank} of the Bohr set, and $\rho$ as the \emph{radius}.
\end{definition}

\begin{example} The arithmetic progression $[N]$ is an uncentred Bohr set of rank $1$, with $S = \{1\}$, $\alpha_1 = (N+1)/2p$ and $\rho = N/2p$.
More generally, any arithmetic progression is an uncentred Bohr set of rank $1$, and conversely.  The intersection of $d$ arithmetic progressions of equal length will be an uncentred Bohr set of rank $d$.  (In fact, in a cyclic group of prime order, this essentially describes all the possible uncentred Bohr sets.)
\end{example}
 
The inverse $U^3$-theorem will only require the centred Bohr sets, but we will need the uncentred Bohr sets in the next section, when we convert the inverse theorem into a Koopman-von Neumann type structure theorem.

Dealing with Bohr sets is slightly technical. One reason for this is that $|B(S,\rho)|$ is not guaranteed to depend particularly smoothly on $\rho$. As discovered by Bourgain \cite{Bou} (see also \cite[Chapter 8]{gt-inverseu3}), such a property \emph{can} be guaranteed for a large supply of $\rho$. To discuss this issue, the following definition is pertinent.

\begin{definition}[Regular Bohr sets]\label{regular-def} Let $S \subseteq \Z/p\Z$ be a set with size $d = |S|$, and suppose that $0 < \rho < 1/2$. A Bohr set $B(S,\rho)$ is said to be \emph{regular} if one has
\[ (1-100d |\kappa|) |B(S,\rho)| \leq |B(S, (1+\kappa)\rho)| \leq (1+100d |\kappa|) |B(S,\rho)|\]
whenever $|\kappa| \leq 1/100d$.
\end{definition}

The \emph{raison d'\^etre} for this definition is a result of Bourgain \cite{Bou} (see also \cite[Lemma 8.2]{gt-inverseu3}) which states that for any $S$ and any $\varepsilon$ there is at least one regular value of $\rho$ in the interval $[\varepsilon,2\varepsilon]$. This will not concern us here though it was important for the proofs in \cite{gt-inverseu3}.  

We move swiftly on to some other concepts which are useful in the discussion of the inverse theorem for the $U^3(\Z/p\Z)$-norm.

\begin{definition}[Linear phase functions]  We say that a function $\phi: \Z/p\Z \to \R/\Z$ is a \emph{globally linear phase function} if we have
\[ \phi(x+h_1+h_2) - \phi(x+h_1) - \phi(x+h_2) + \phi(x) = 0\]
for all $x,h_1,h_2 \in \Z/p\Z$.
\end{definition}

\begin{example}
Because $p$ is prime, it is easy to see that a function $\phi: \Z/p\Z \to \R/\Z$ is globally linear if and only if it takes the form $\phi(x) = \xi x/p + \alpha$ for some $\xi \in \Z/p\Z$ and $\alpha \in \R/\Z$.
\end{example}

\begin{definition}[Quadratic phase functions]  Let $B \subset \Z/p\Z$.  We say that a function $\phi: B \to \R/\Z$ is a \emph{locally quadratic phase function} on $B$ if we have
\begin{equation}\label{quaddef}
\begin{split}
 \phi(x+h_1+h_2+h_3) &- \phi(x+h_1+h_2) - \phi(x+h_2+h_3) - \phi(x+h_1+h_3) \\
 &+ \phi(x+h_1) + \phi(x+h_2) + \phi(x+h_3) - \phi(x) = 0
\end{split}
\end{equation}
whenever $x$, $x+h_1$, $x+h_2$, $x+h_3$, $x+h_1+h_2$, $x+h_1+h_3$, $x+h_2+h_3$, and $x + h_1 + h_2 + h_3$ all lie in $B$.
\end{definition}

\begin{example}  Every globally linear phase function is locally quadratic.
If $\alpha,\beta,\gamma$ are real numbers, and $N$ is an integer between $p/8$ and $p/4$,
then the function $\phi(n) = \alpha n^2 + \beta n + \gamma \md{1}$ is a locally quadratic phase function on $[N]$.
\end{example}

\begin{remark}
There are also notions of locally linear phase functions, and globally quadratic ones, but we will not need them here.
\end{remark}

We are now ready to state the inverse theorem for the $U^3(\Z/p\Z)$-norm in the form that we shall need it. We write $e(x) := e^{2\pi i x}$ as usual.

\begin{theorem}[Inverse $U^3(\Z/p\Z)$ theorem] \label{invert-u3}
Let $f: \Z/p\Z \to \C$ be a $1$-bounded function 
such that $\|f\|_{U^3(\Z/p\Z)} \geq \eta$ for some $\eta \in (0,1)$.  Then there exists a regular Bohr set $B := B(S,\rho)$
with $|S| \ll \eta^{-C}$ and $\rho \gg \eta^C$, and a locally quadratic phase function $\phi_y: y + B \to \R/\Z$
on $y + B$ for every $y \in \Z/p\Z$, such that
\begin{equation}\label{efphi-general}
\E_{y \in \Z/p\Z} |\E_{t \in y+B} f(t) e(-\phi_y(t))| \gg \eta^C.
\end{equation}
\end{theorem}

This is \cite[Theorem 2.7]{gt-inverseu3}, where in fact the explicit value of $C = 2^{24}$ was attained.

\begin{remark} It is unfortunately necessary to deal with locally quadratic phase functions rather than the more intuitively natural
globally quadratic phase functions; see \cite{gowers-4-aps,green-icm,gt-inverseu3} for further discussion of issues of this type, or
the paper \cite{furstenberg-von-neumann} for a rather different perspective on the same phenomenon.
\end{remark}

\section{Linear and quadratic factors, and a quadratic Koopman-von Neumann theorem}\label{kvn-sec}

As in \cite{green-tao-szem}, we now use an ``energy increment argument'' to convert our inverse theorem to a quadratic structure theorem of Koopman-von Neumann type, 
inspired by some ideas from ergodic theory. Part I of the series \cite{green-tao-szem} or the lecture notes \cite{green-montreal} may be consulted for further discussion, and \cite{tao:icm} gives a more general discussion of structure theorems and inverse theorems.  We first need some more notation.

\begin{definition}[Factors]  Let $W$ be any non-empty finite set.
Define a \emph{factor} (or \emph{$\sigma$-algebra})
in $W$ to be a collection $\B$ of subsets of $W$ which are closed under union, intersection, and complement, and which contains $\emptyset$ and $W$. Define an \emph{atom} of $\B$ to be a minimal non-empty subset of $W$; these partition $W$, and indeed in this finitary setting a factor may be thought of simply as a partition of $W$.  If $\B$, $\B'$ are factors in $W$ with $\B \subseteq \B'$ we say that $\B'$ \emph{extends} $\B$.  More generally, if $\B, \B'$ are factors in $W$ we let $\B \vee \B'$ be the smallest common extension, so that the atoms of $\B \vee \B'$ are the intersections of atoms of $\B$ and atoms of $\B'$.  If $\B$ is a factor in $W$ and $W'$ is a subset of $W$, we define the \emph{restriction}
$\B|_{W'}$ of $\B$ to $W'$ to be the factor of $W'$ formed by intersecting all the sets in $\B$ with $W'$.
If $f: W \to \C$, we let $\E(f|\B): W \to \C$ denote the conditional expectation
$$ \E(f|\B)(x) := \E(f|\B(x)) \hbox{ for all } x \in W,$$
where $\B(x)$ is the unique atom in $\B$ that contains $x$.  Equivalently, $\E(f|\B)$ is the orthogonal projection to the space $\B$-measurable functions in the Hilbert space
$L^2(\Z/p\Z)$.  
\end{definition}

We will focus our attention on very structured factors, namely linear and quadratic factors, which are generated from globally linear
and locally quadratic phase functions respectively.  The notation here is inspired by the finite field analogues in 
\cite{green-tao-szem} but with one new parameter, a ``resolution'' $K$, which is needed as a substitute for the small torsion that one 
enjoys in the finite field geometry setting. We first need to describe how to convert a phase function into a factor.

\begin{definition}\label{def-52} Call a phase function \emph{irrational} if it only takes irrational values.  If $\phi: W \to \R/\Z$ is an irrational 
phase function on a finite nonempty set $W$ and $K \geq 1$ is an integer, we define $\B_{\phi,K}$ to be the factor in $W$ whose atoms are the sets
$\{ x \in \Z/p\Z: \| \phi(x) - j/K \|_{\R/\Z} < 1/2K \}$ for $j=0,1,\ldots,K-1$.
\end{definition}

\begin{remark} The assumption of irrationality is a minor technicality, used in order to avoid having to deal with the borderline case when $\| \phi(x) - j/K \|_{\R/\Z}$ is exactly equal to $1/2K$; in practice we shall be able to use perturbation arguments to work purely with irrational phase functions.
\end{remark}

\begin{definition}[Linear factors]  A \emph{linear factor of complexity at most $d$ and resolution $K$} is any factor $\B$ in $\Z/p\Z$
of the form $\B = \B_{\phi_1,K} \vee \ldots \vee \B_{\phi_{d'},K}$, where $0 \leq d' \leq d$ and $\phi_1, \ldots, \phi_{d'}: \Z/p\Z \to \R/\Z$ are irrational globally linear phase functions.
\end{definition}

\begin{remark} From the definitions we see that if $\B$ is a linear factor of complexity at most $d$ and resolution $K$, then $\B$ has at most $K^d$ atoms, each of which is an uncentred Bohr set of rank at most $d$ and radius $1/2K$.  Also, if $\B'$ is another linear factor of complexity at most $d'$ and resolution $K$, then clearly
$\B \vee \B'$ is a linear factor of complexity at most $d+d'$ and resolution $K$.
\end{remark}

\begin{definition}[Quadratic factors]  Let $B$ be an uncentred Bohr set.  A \emph{pure quadratic factor of complexity at most $d$ and resolution $K$}
in $B$ is any factor $\B$ in $B$ of the form $\B = \B_{\phi_1,K} \vee \ldots \vee \B_{\phi_{d'},K}$, where $0 \leq d' \leq d$ and $\phi_1, \ldots, \phi_d: B \to \R/\Z$ are irrational locally quadratic phase functions on $B$.  A \emph{quadratic factor of complexity at most $(d_1,d_2)$ and resolution $K$} is any pair $(\B_1,\B_2)$ of factors in $W$, where $\B_1$ is a linear factor of complexity at most $d_1$ and resolution at most $K$, and $\B_2$ is an extension of $\B_1$, whose restriction to any atom $B$ of $\B_1$ is a pure quadratic factor on $B$ of complexity at most $d_2$ and resolution at most $K$.  
We say that one quadratic factor $(\B'_1,\B'_2)$ is a \emph{quadratic extension} of another $(\B_1,\B_2)$ if 
$\B_1 \subseteq \B'_1$ and $\B_2 \subseteq \B'_2$.
\end{definition}

\begin{remark}
Observe that if $(\B_1,\B_2)$ and $(\B'_1,\B'_2)$ are quadratic factors of resolution $K$ and complexity at most $(d_1,d_2)$ and $(d'_1,d'_2)$ respectively,
then their common extension $(\B_1 \vee \B'_1, \B_2 \vee \B'_2)$ is a quadratic factor of complexity at most $(d_1+d'_1,d_2+d'_2)$; this is ultimately
because the restriction of a locally quadratic phase function to a smaller set remains locally quadratic.\end{remark}
Our next task is to rephrase the inverse theorem, Theorem \ref{invert-u3}, in terms of quadratic factors. At heart this is really nothing more than an averaging argument, though due to ``edge effects'' it is somewhat tedious to write down rigorously.

\begin{theorem}[Inverse theorem for $U^3(\Z/p\Z)$, again]\label{invert-u3-sigma}
Let $f: \Z/p\Z \to \C$ be a $1$-bounded function such that $\|f\|_{U^3(\Z/p\Z)} \geq \eta$ for some $\eta \in (0,1)$.  Suppose also
that $K$ is an integer such that $K \geq C \eta^{-C}$ for some sufficiently large constant $C > 0$.  Then there exists a quadratic factor $(\B_1,\B_2)$
in $\Z/p\Z$ of complexity at most $(O( \eta^{-C} ), 1)$ and resolution $K$ such that
\begin{equation}\label{fib}
\| \E(f|\B_2) \|_{L^1(\Z/p\Z)} \gg \eta^C.
\end{equation}
\end{theorem}

\begin{proof}  Let $S, \rho$ be as in Theorem \ref{invert-u3}.  
Let $\alpha = (\alpha_\xi)_{\xi \in S}$ be a point on the torus $(\R/\Z)^S$ with irrational coefficients (one could chose it randomly, if desired).  
We then define $\B_1$ to be the $\sigma$-algebra whose atoms are of the form
$$ \{ x \in \Z/p\Z: \| x \xi/p - \alpha_\xi - j_\xi / K \|_{\R/\Z} < 1/2K \hbox{ for all } \xi \in S \}$$
where for each $\xi \in S$, $j_\xi$ is an integer between $0$ and $K-1$.  One easily verifies that $\B_1$ is an irrational linear factor of complexity $|S|$ and resolution $K$, defined by linear phases $\phi_{\xi}(x) := x\xi/p - \alpha_{\xi}$, $\xi \in S$.  For each $x \in\Z/p\Z$, let $F(x)$ be the quantity
\[ F(x) := \sup_{\phi \in \Phi(x)} |\E_{t \in \B_1(x)}( f(t) e(-\phi(t)) )|,\]
where $\B_1(x)$ is the atom of $\B_1$ that contains $x$ and $\Phi(x)$ is the collection of all locally quadratic phase functions $\phi: \B_1(x) \to \R/\Z$. Thus $F$ measures the maximum correlation of $f$ with a quadratic phase on the atom $\B_1(x)$.  We claim that it suffices to show that 
\begin{equation}\label{fl1}
\| F \|_{L^1(\Z/p\Z)} \gg \eta^C.
\end{equation}
Suppose that this has been established. Written out in full, it becomes the statement that
\[ \E_{x \in \Z/p\Z} \big| \E_{t \in \B_1(x)} f(t) e(-\phi_{\B_1(x)}(t))  \big| \gg \eta^C\] for an appropriate choice of $\phi_{\B_1(x)} \in \Phi(x)$.
Modulating each phase by a complex number $e(\theta_x)$ we may move the modulus signs to the outside, obtaining
\[ \big|\E_{x \in \Z/p\Z} \E_{t \in \B_1(x)} f(t) e(-\phi_{\B_1(x)}(t))  \big|\gg \eta^C.\]
Since the two averaging operations are equivalent to the single averaging $\E_{x \in \Z/p\Z}$ this becomes
\begin{equation}\label{use-soon-1} |\E_{x \in\Z/p\Z} f(x) e(-\phi_{\B_1(x)}(x)) | \gg \eta^C.\end{equation}
By perturbing each of the $\phi_{B_1}$ infinitesimally we may assume that the $\phi_{B_1}$ are all irrational.  If we then let $\B_2$ be the extension of $\B_1$ whose restriction to each atom $B_1$ of $\B_1$ is given by $\B_{\phi_{B_1},K}$, then $(\B_1,\B_2)$ is a quadratic factor of complexity
at most $(|S|, 1)$ and resolution $K$, and we have
\[ e(-\phi_{\B_1(x)}(x)) = \E( e(-\phi_{\B_1(x)}) | \B_2 )(x)+ O(1/K)\]
for all $x$. It is important to note here that $\phi_{\B_1(x)}$ depends only on the atom $\B_1(x)$ and not otherwise on $x$ itself.

It follows from this and \eqref{use-soon-1} that if  $K \geq C \eta^{-C}$ for sufficiently large $C$ then
\[ \big|\langle f, \E( e(\phi_{\B_1(x)})|\B_2)\rangle\big| \gg \eta^C,\]
where we have written $\langle g_1,g_2 \rangle := \E_{x \in \Z/p\Z} g_1(x) \overline{g_2(x)}$.
The conditional expectation operator $g \mapsto \E(g | \B_2)$ is self-adjoint with respect to this inner product, and hence this implies that
\[ \big| \langle \E(f | \B_2), e(\phi_{\B_1(x)}) \rangle \big| \gg \eta^C.\]
The desired bound \eqref{fib} is now a consequence of the triangle inequality in the form $|\langle g_1,g_2\rangle| \leq \Vert g_1 \Vert_1 \Vert g_2 \Vert_{\infty}$.

It remains, then, to establish \eqref{fl1}. It is now time to exploit the estimate \eqref{efphi-general}, which we urge the reader to recall now.  For any fixed $y \in \Z/p\Z$, we write $\Omega_y$ for the union of those atoms of $\B_1$ which only partially intersect $y + B$ (thus they are neither contained in $y + B(S,\rho)$ nor outside of it).  We have
\begin{align}\nonumber
|\E_{t \in y+B} f(t) e(-\phi_y(t))| &\leq \sum_{B_1: B_1 \subset y+B} \frac{|B_1|}{|y+B|} |\E_{t \in B_1} f(t) e(-\phi_y(t))| + \P_{y+B}(\Omega_y) \\ \nonumber
&\leq \sum_{B_1: B_1 \subset y+B} \frac{|B_1|}{|y+B|} \E_{B_1} F + \P_{y+B}(\Omega_y) \\
&\leq \E_{y+B} F + \P_{y+B}(\Omega_y)\label{previous-bound}
\end{align}
(note that $F$ is constant on each atom $B_1$).  
On the other hand, since $B = B(S,\rho)$, one easily verifies that
\[ \Omega_y \subseteq y + \big(B(S,\rho + \frac{1}{K}) \backslash B(S, \rho - \frac{1}{K})\big).\]
since $B(S,\rho)$ is regular, we conclude that
\[ \P_{y+B}(\Omega_y) \ll \frac{d}{\rho K}.\]
Inserting this into \eqref{previous-bound}, taking expectations, and using \eqref{efphi-general}, we conclude that
\[ \| F \|_{L^1(\Z/p\Z)} + O( \frac{d}{\rho K} ) \gg \eta^C.\]
Now $d \ll \eta^{-C}$ and $\rho \gg \eta^C$. Thus by taking $K \geq C\eta^{C'}$ for large enough $C'$ we obtain the claim that $\Vert F \Vert_{L^1(\Z/p\Z)} \gg \eta^C$.\end{proof}

Let $\B_{\triv}$ denote the rather trivial factor generated by the two atoms $[N]$ and $(\Z/p\Z) \setminus [N]$. With $f$ as in Theorem \ref{dens-inc} and this new notation we have $\delta 1_{[N]} = \E(f|\B_\triv)$. Our next task is to iterate Theorem \ref{invert-u3-sigma} 
via an energy increment argument to obtain the following structural result of ``Koopman-von Neumann'' type. The blueprint for arguments of this type is Szemer\'edi's proof of his \emph{regularity lemma} in graph theory \cite{szem:reg}. For other examples in additive combinatorics the reader might consult any of \cite{green-fin-field,green-montreal,green-tao-primes,green-tao-szem,tao-survey,taovu-book}.

\begin{theorem}[Quadratic Koopman-von Neumann theorem]\label{kvn}  
Let $f: \Z/p\Z \to [-1,1]$ be a $1$-bounded function, and let $\eta > 0$.  
Suppose also that $K$ is an integer such that $K \geq C \eta^{-C}$ for some sufficiently large constant $C > 0$.  Then there exists a
quadratic factor $(\B_1,\B_2)$ in $\Z/p\Z$ of complexity at most $(O( \eta^{-C} ), O( \eta^{-C} ))$
and resolution $K$ such that
\begin{equation}\label{efb-b1}
\|f-\E(f|\B_2 \vee \B_\triv)\|_{U^3(\Z/p\Z)} \leq \eta.
\end{equation}
\end{theorem}

\begin{proof} We run the following algorithm:
\begin{itemize}
\item Step 0: Initialize $\B_1 = \B_2 = \{ \emptyset, \Z/p\Z\}$. Thus $(\B_1,\B_2)$ is a quadratic factor with complexity $(0,0)$ and resolution $K$.
\item Step 1: If \eqref{efb-b1} holds, then \textsc{stop}.  Otherwise, apply Theorem \ref{invert-u3-sigma} with $f$ replaced\footnote{This function is bounded pointwise by 2. It is clear that a trivial rescaling of Theorem \ref{invert-u3-sigma} applies to such functions.} by $f - \E(f|\B_2 \vee \B_\triv)$ to obtain a quadratic factor $(\B'_1, \B'_2)$ of complexity at most $(O( \eta^{-C} ), 1)$ and resolution $K$ such that
\begin{equation}\label{l1proj}
\| \E( f - \E(f|\B_2 \vee \B_\triv) | \B'_2 ) \|_{L^1(\Z/p\Z)} \gg \eta^C.
\end{equation}
\item Step 2: Replace $(\B_1,\B_2)$ with $(\B'_1, \B_1 \vee \B'_2)$ (thus increasing the complexity of $(\B_1,\B_2)$ by at most $(O( \eta^{-C} ), 1)$),
and return to Step 1.
\end{itemize}

Observe from \eqref{l1proj} and Cauchy-Schwarz that
\[ \| \E( f - \E(f|\B_2 \vee \B_\triv) | \B'_2 ) \|_{L^2(\Z/p\Z)} \gg \eta^C,\]
and hence 
$$ \| \E( f - \E(f|\B_2 \vee \B_\triv) | \B_2 \vee \B'_2 \vee \B_{\triv} ) \|_{L^2(\Z/p\Z)} \gg \eta^C.$$
By Pythagoras' theorem we conclude that
\[ \| \E(f | \B_2 \vee \B'_2 \vee \B_\triv ) \|_{L^2(\Z/p\Z)}^2 - \| \E(f | \B_2 \vee \B_\triv ) \|_{L^2(\Z/p\Z)}^2
\gg \eta^C.\]
It follows that every time we perform Step 2, the energy $\| \E(f | \B_2 \vee \B_\triv ) \|_{L^2(\Z/p\Z)}^2$ increments by at least $\gg \eta^C$.
Since the energy is clearly bounded between $0$ and $1$, the algorithm can only run for at most $O( \eta^{-C})$ iterations, and the claim easily follows.
\end{proof}

If we apply this theorem (with $\eta := c \delta^4$ for some small $c$) and Lemma \ref{u3-lambda4} to the situation in Theorem \ref{dens-inc} 
we obtain the following corollary.

\begin{corollary}[Anomalous AP4 count on a quadratic factor]\label{dens-cor}
Let the assumptions be as in Theorem \ref{dens-inc}.  
Then there exists a quadratic factor $(\B_1,\B_2)$ in $\Z/p\Z$ of complexity at most $(O(\delta^{-C}), O(\delta^{-C}))$
and resolution $O( \delta^{-C})$
such that the function $g := \E(f|\B_2 \vee \B_\triv)$ obeys
\begin{equation}\label{gap}
 |\Lambda(g,g,g,g) - \Lambda(\delta 1_{[N]}, \delta 1_{[N]}, \delta 1_{[N]}, \delta 1_{[N]})| \gg \delta^4.
 \end{equation}
\end{corollary}

Thus we have replaced the original function $f$ by the more structured function $g$.  Note that $\delta 1_{[N]} = \E(g|\B_{\triv})$. From this it is not hard to obtain, under the assumption that $f$ has an anomalous count of 4-term progressions, a substantial density increment for $f$ 
on a quadratic Bohr set.

\begin{corollary}[Density increment on quadratic Bohr set]\label{dens-cor2}
Let the assumptions be as in Theorem \ref{dens-inc}.  
Then there exists a quadratic factor $(\B_1,\B_2)$ in $\Z/p\Z$ of complexity at most $(O(\delta^{-C}), O(\delta^{-C}))$
and resolution $O( \delta^{-C}) $, and an atom $B_2$ of $\B_2 \vee \B_\triv$ of density $\P_{\Z/p\Z}(B_2) \gg \exp( - O( \delta^{-C} ) )$
and contained in $[N]$ such that
\[ \E_{B_2}(f) \geq (1 + c) \delta\]
for some absolute contant $c > 0$.
\end{corollary}

\begin{proof}  Let $(\B_1,\B_2)$ and $g$ be as in Corollary \ref{dens-cor}. The facts that $[N]$ is measurable in $\B_2 \vee \B_\triv$ and that $f$ is supported on $[N]$ guarantee that $g$ is also supported on $[N]$.  Let $\Omega$ denote the set where $g \geq (1+c)\delta$, where $c > 0$ is a small constant to be chosen later, and
let $g' := (1-1_\Omega) g$. From Lemma \ref{l1-lambda4} we have
\[ |\Lambda(g,g,g,g) - \Lambda(g',g',g',g')| \leq 4 \P_{\Z/p\Z}(\Omega)\]
and
\[ |\Lambda(g',g',g',g') - \Lambda(\delta 1_{[N]}, \delta 1_{[N]}, \delta 1_{[N]}, \delta 1_{[N]})| \leq 8 \delta^3
\| g' - \delta 1_{[N]} \|_{L^1(\Z/p\Z)}.\]
Furthermore we evidently have \[ \| g' - \delta 1_{[N]} \|_{L^1(\Z/p\Z)} \leq \| g - \delta 1_{[N]} \|_{L^1(\Z/p\Z)} + \P_{\Z/p\Z}(\Omega).\]  Combining
these three estimates together with \eqref{gap} we obtain
\[ \delta^3 \| g - \delta 1_{[N]} \|_{L^1(\Z/p\Z)} + \P_{\Z/p\Z}(\Omega) \gg \delta^4.\]
Now observe that the positive part $(g - \delta 1_{[N]})_+$ of $g - \delta 1_{[N]}$ can only exceed $c\delta$ on $\Omega$, and hence has a total $L^1$ norm of
at most $c\delta + \P_{\Z/p\Z}(\Omega)$.  Since $g - \delta 1_{[N]}$ also has mean zero, we conclude that
\[\| g - \delta 1_{[N]} \|_{L^1(\Z/p\Z)} = 2\| (g - \delta 1_{[N]})_+ \|_{L^1(\Z/p\Z)} \ll c\delta + \P_{\Z/p\Z}(\Omega).\]
If $c$ is chosen small enough, we deduce that
\[ \P_{\Z/p\Z}(\Omega) \gg \delta^4.\]
Now $\B_2 \vee \B_\triv$ has complexity and resolution $O(\delta^{-C})$, and hence contains at most $\exp( O( \delta^{-C} ) )$ atoms. By the pigeonhole principle we can therefore find an atom $B_2$ of $\B_2 \vee \B_\triv$ which is contained in $\Omega$ and which has 
$\P_{\Z/p\Z}(B_2) \gg \exp( - O( \delta^{-C} ) )$.  By construction we have $\E_{B_2}(f) \geq (1+c)\delta$, and the claim follows.\end{proof}

Our sole remaining task is to take this density increment for $f$ on a ``quadratic Bohr set'' and use it to obtain a similar density increment for $f$ on an arithmetic progression (Theorem \ref{dens-inc}). This we do by splitting the quadratic Bohr set into a union of progressions, a process we call \emph{linearisation}.

\section{Linearisation of quadratic Bohr sets}

We will decompose a quadratic Bohr set into a union of progressions. Our method for doing this does not naturally output progressions of equal sizes, and the following simple variant of the pigeonhole principle is designed to ensure that there is at least one progression which is quite long and on which $f$ has a substantial density increment.

\begin{lemma}[Pigeonhole principle]\label{dens-split}  
Let $B$ be a non-empty set, and let $B = A_1 \cup \ldots \cup A_m$ be a partition of $B$ into $m$ disjoint sets.  Let $f: B \to \R^+$ be a $1$-bounded
nonnegative function.  Then for any $\eps > 0$, there exists $i \in \{1,\dots,m\}$ such that $\P_B(A_i) > \eps/m$ and
$$ \E_{A_i}(f) \geq \E_B(f) - \eps.$$
\end{lemma}

\begin{proof}  Let $\Omega$ be the union of all the $A_i$ for which $\P_B(A_i) \leq \eps/m$. We obviously have $\P_B(\Omega) \leq \eps$.  From Bayes' identity and the fact that $0 \leq f \leq 1$ we have
\[ \E_B(f) = \P_B(\Omega) \E_\Omega(f) + (1-\P_B(\Omega)) \E_{B \backslash \Omega}(f) \leq \P_B(\Omega) + \E_{B \backslash \Omega}(f), \] and it follows that
\[ \E_{B \backslash \Omega}(f) \geq \E_B(f) - \eps.\]
Partitioning $B \backslash \Omega$ into its constituent sets $A_i$, the claim follows from the usual pigeonhole principle.
\end{proof}

The next result provides the splitting of a quadratic Bohr sets into progressions, and is the main result of the section.

\begin{proposition}[Linearisation of quadratic Bohr sets]\label{linquad} Let $(\B_1,\B_2)$ in $\Z/p\Z$ of complexity at most $(d_1,d_2)$
and some resolution $K$, and let $B_2$ be an atom of $\B_2$.  Then one can partition $B_2 \cap [N]$ as the union of 
$\ll d_2^{O(d_2)} N^{1 - c/(d_1+1)(d_2+1)^3}$ disjoint arithmetic progressions in $\Z/p\Z$.  
\end{proposition}

\emph{Proof of Theorem \ref{dens-inc} assuming Proposition \ref{linquad}.} Suppose that $f$ satisfies the conditions of Theorem \ref{dens-inc}. By Corollary \ref{dens-cor2} we know that there is a quadratic factor $(\B_1,\B_2)$ in $\Z/p\Z$ of complexity at most $(O(\delta^{-C}), O(\delta^{-C}))$ and resolution $O(\delta^{-C})$, and an atom $B_2 \subseteq [N]$ of $\B_2 \vee \B_\triv$, having density at least $\exp(-O(\delta^{-C}))$ in $\Z/p\Z$, on which the average of $f$ is at least $(1 + c_0)\delta$. Using Proposition \ref{linquad} we may write $B_2$ as the union of $\exp(O(\delta^{-C})) N^{1 - c\delta^C}$ progressions. Taking $\eps := c_0\delta/2$ in Lemma \ref{dens-split}, we obtain a progression of length at least $\exp(-O(\delta^{-C})) N^{c\delta^C}$ on which $f$ has average at least $(1 + \frac{1}{2}c_0)\delta$. To complete the proof of the theorem, we need to make sure that this length is in fact $\gg N^{c'\delta^{C'}}$ for absolute constants $c',C' > 0$. This may be ensured by taking the absolute constant in the condition \eqref{plarge} to be sufficiently large.\endproof 

It remains to prove Proposition \ref{linquad}.  We first deal with the linear component of the factor $(\B_1,\B_2)$.  Since $\B_2$ extends $\B_1$, there is a unique atom $B_1$ in $\B_1$ which contains $B_2$.

\begin{proposition}[Linearisation of linear Bohr sets]\label{linlin} Let $\B_1$ be a linear factor of complexity $d_1$ and resolution $K$. Let $B_1$ be an atom in $\B_1$. Then one can partition $B_1 \cap [N]$ as the union of $\ll 2^{d_1}N^{1 - 1/(d_1+1)}$ arithmetic progressions.
\end{proposition}

\begin{proof}  We can write $B_1$ as an uncentred Bohr set $B_\alpha(S, 1/2K)$, where $|S| \leq d_1$ and $\alpha \in (\R/\Z)^S$.  Using the Kronecker approximation theorem (Proposition \ref{kat}) we can find a non-zero $r \in \Z/p\Z$ such that
\[ \| \xi r / p \|_{\R/\Z} \ll N^{-1/(d_1+1)} \]
for all $\xi \in S \cup \{1\}$.  If we then partition $\Z/p\Z$ into $O( N^{1 - 1/(d_1+1)} )$ arithmetic progressions of common difference $r$ and length
$O( N^{1/(d_1+1)} )$, we see that the intersection of each of these progressions with $B_1 \cap [N]$ will be the union of no more than $2^{d_1}$ smaller arithmetic
progressions, also of step $r$, and the claim follows.
\end{proof}

This last proposition improves our situation considerably, since it is much easier to understand quadratic phases on a progression than it is quadratic phases on a Bohr set.

\begin{proposition}[Linearisation of pure quadratic Bohr sets]\label{linpure}  Let $P$ be an arithmetic progression in $\Z/p\Z$, and let $\phi_1, \ldots, \phi_d: P \to \R/\Z$ be locally quadratic irrational phase functions on $P$. Consider the factor $\B_{\phi_1,K} \vee \ldots \vee \B_{\phi_d,K}$ of resolution $K$ defined by these phase functions \textup{(}cf. Definition \ref{def-52}\textup{)}. Then for any resolution $K$, every atom $B_2 \subseteq P$ 
of  can be partitioned as the union of $\ll d^{O(d)} |P|^{1 - c/(1+d)^3}$ disjoint arithmetic progressions.
\end{proposition}

Let us now see why Proposition \ref{linlin} and Proposition \ref{linpure} together imply Proposition \ref{linquad}.  First let us take all the progressions arising from Proposition \ref{linlin} which are rather small, say having length $O( N^{1/2(d_1+1)})$.  We can partition these progressions in the most trivial way into singletons, ending up with at most $O( N^{1 - 1/2(d_1+1)} )$ single-element
progressions in this way.  As for each longer progression $P$ in $B_1 \cap [N]$, of length $\gg N^{1/2(d_1+1)}$, we apply Proposition \ref{linpure} to see that $P \cap B_2$ is the union of $\ll d_2^{O(d_2)} |P|^{1 - c/(1+d_2)^3} \ll d_2^{O(d_2)} N^{-c/(1+d_1)(1+d_2)^3} |P|$ disjoint arithmetic progressions.  Assembling all of these progressions together as $P$ varies, we obtain $O( d_2^{O(d_2)} N^{1-c/(1+d_1)(1+d_2)^3})$ disjoint progressions in total, and Proposition \ref{linquad} follows.

It remains to prove Proposition \ref{linpure}. A result of this type, in which there is just a single quadratic phase, may be found in \cite{gowers-4-aps}. Here, however, we are dealing with $d$ quadratics rather than just one, and will have to take a little care to make sure that our exponents depend only polynomially on $d$ rather than exponentially. Because of this, we cannot, for example, simply iterate the analogous single-quadratic results from
\cite{gowers-4-aps}.  As a first step we may apply an affine linear transformation to $P$ and assume that $P = [1,M]$ for
some $M$, $1 \leq M \leq p$. We can also take $d \geq 1$ since the $d = 0$ case is trivial. It is easy to see, straight from the definition of a quadratic phase, that each $\phi_j: [1,M] \to \R/\Z$ takes the form
\[
\phi_j(n) = \alpha_j n^2 + \beta_j n + \gamma_j\]
for some $\alpha_j, \beta_j, \gamma_j \in \R/\Z$.  The set $B_2$ thus takes the form
$$ \{ n \in [1,M]: \| \alpha_j n^2 + \beta_j n + \gamma'_j \|_{\R/\Z} < 1/2K \hbox{ for } j = 1,\dots,d \}$$
where $\gamma'_j$ is some other element of $\R/\Z$.  Our objective is to partition this set into $\ll d^{O(d)} M^{1-c/d^3}$ disjoint arithmetic
progressions.

The first step, as in \cite{gowers-4-aps}, is to find a scale $r$ for which the effects of the quadratic components $\alpha_j n^2$ of each 
phase are locally negligible.  Applying Proposition \ref{sqj} we can locate an integer $r$, $1 \leq r \leq \sqrt{M}$, such that 
\begin{equation}\label{alphaj-quadratic}
\| \alpha_j r^2 \|_{\R/\Z} \ll d M^{-c_0/d^2} 
\end{equation}
whenever $1 \leq j \leq d$, where $c_0 > 0$ is an absolute constant.
Now we can partition $[1,M]$ into at most $M^{1-c_0/4d^2}$ arithmetic progressions of step $r$ and lengths $\sim M^{c_0/4d^2}$ (that is, bounded above and below by absolute constants times this). It will suffice
to show that, for each such arithmetic progression $P$, the set $P \cap B_2$ can be partitioned into $\ll d^{O(d)} |P|^{1-1/2d}$ arithmetic progressions.

Let us fix one of these progressions $P = a, a+r, \ldots, a+(k-1)r$, where $k \sim M^{c_0/4d^2}$. From \eqref{alphaj-quadratic} we have
\begin{equation}\label{alphaj-quadratic2}
\| \alpha_j r^2 \|_{\R/\Z} \ll d k^{-4} \hbox{ for } j = 1,\dots,d. 
\end{equation}
The set $P \cap B_2$ can be written as
$$ \{ a+ir: \| i^2 \alpha_j r^2 + \beta_{j,P} i + \gamma_{j,P} \|_{\R/\Z} < 1/2K \hbox{ for } j = 1,\dots,d \}$$
where $\beta_{j,P}, \gamma_{j,P}$ are some real numbers depending on $j$ and $P$.  Now we use Kronecker's theorem (Proposition \ref{kat}) to find a positive integer $s \leq \sqrt{k}$ such that
\begin{equation}\label{betajp-linear}
 \| \beta_{j,P} s \|_{\R/\Z} \ll k^{-1/2d}
 \end{equation}
for all $j \in \{1,\dots,d\}$.  We now partition $P$ into $\ll k^{1-1/2d}$ arithmetic progressions of step $rs$ and length $\ll k^{1/2d}$.  Consider a single such progression $Q$.  It can be written as
\[ Q = \{ a + (b + ts)r: 1 \leq t \leq T \}\]
for some $b \leq k$ and some $T \ll k^{1/2d}$, and its intersection with $B_2$ can be written as
\[ \{ a + (b + ts)r: 1 \leq t \leq T; 
\| (b+ts)^2 \alpha_j r^2 + \beta_{j,P} (b+ts) + \gamma_{j,P} \|_{\R/\Z} < 1/2K \hbox{ for } j = 1,\dots,d\}.\]
The expression $(b+ts)^2 \alpha_j r^2 + \beta_{j,P} (b+ts) + \gamma_{j,P}$ can be rewritten (modulo $1$) as
$$ t^2 s^2 \{ \alpha_j r^2 \} + t (2bs \{ \alpha_j r^2 \} + \{ \beta_{j,P} s \}) + c_{j,P,Q}$$
where $\{x\} \in (-1/2,1/2]$ is the difference between $x$ and the nearest integer to $x$, and $c_{j,P,Q}$ is some real number.  Observe from \eqref{alphaj-quadratic2}, \eqref{betajp-linear} and the bounds $s \leq \sqrt{k}$, $b \leq k$, $T \ll k^{1/2d}$ that
the coefficients of $t^2$ and $t$ in this quadratic polynomial are $O(d/T^2)$ and $O(d/T)$ respectively.  Thus, for each fixed $j$, the set of values $t$ for which this expression has an $\R/\Z$ norm less than $1/2K$ is the union of $O(d)$ intervals (arithmetic progressions of step $1$). 
This means that $Q \cap B_2$ is the union of at most $O(d)^d \ll d^{O(d)}$ intervals, and thus $P \cap B_2$ can be partitioned into
$\ll d^{O(d)} k^{1-1/2d}$ progressions as desired.  This concludes the proof of Proposition \ref{linpure} and hence, by earlier reductions, that of our main theorem.\endproof

\appendix
\section{Simultaneous quadratic recurrence}

We recall the well-known Kronecker approximation theorem:

\begin{proposition}[Kronecker approximation theorem]\label{kat}  Let $\alpha_1, \ldots, \alpha_d$ be real numbers, and let $N \geq 1$ be an integer.  Then there exists an integer $n$, $1 \leq n \leq N$, such that
\begin{equation}\label{nj-linear}
\| n \alpha_j \|_{\R/\Z} \ll N^{-1/d} \hbox{ for } j = 1,\dots,d.
\end{equation}
\end{proposition}

This is easily deduced from the pigeonhole principle, partitioning the torus $(\R/\Z)^d$ into fewer than $N$ regions of diameter $O(N^{1/d})$ each, and considering the orbit of $(n\alpha_1,\ldots,n\alpha_d)$. The objective of this appendix is to prove the following quadratic analogue of the above theorem, due to Schmidt \cite{schmidt}.

\begin{proposition}[Simultaneous quadratic recurrence]\label{sqj}  Let $\alpha_1, \ldots, \alpha_d$ be real numbers, and let $N \geq 1$ be an integer.  Then there exists an integer $1 \leq n \leq N$ such that
\begin{equation}\label{nj}
\| n^2 \alpha_j \|_{\R/\Z} \ll d N^{-c/d^2} \hbox{ for } j = 1,\dots,d.
\end{equation}
Here $c > 0$ is an absolute constant.
\end{proposition}

In actual fact Schmidt shows that one may satisfy
\[ \|n^2 \alpha_j\|_{\R/\Z} \ll_{d,\eps} N^{-1/(d^2 + d) + \eps}.\]
The exponent here is of course more precise than the one we quote, but it is of critical importance for our work that we have some understanding of the dependence on $d$ of the implied constant in the $\ll_{d,\eps}$. We could allow it to be (say) $e^{O(d^C)}$, but not much worse. Schmidt's argument is explicit and effective enough that such bounds can probably be extracted with some effort from \cite{schmidt}; but for the convenience of the reader we shall instead provide a complete and self-contained proof of Proposition \ref{sqj} in this appendix.  Note that we only require an exponent of shape $N^{-c/d^C}$ in \eqref{nj}, which is somewhat weaker than what \cite{schmidt} gives, but we do not know of a way to obtain such an exponent which does not follow Schmidt's argument. An exponent $N^{-1/C^d}$ may be obtained by the simpler device of iteratively applying the case $d = 1$ of Proposition \ref{sqj} (see \cite{green-square} for details), but this does not suffice for our purposes here.

Let us begin by sketching some features of Schmidt's argument. Suppose one wishes to find an $n \leq N$ such that $\Vert n \alpha_j \Vert_{\R/\Z} \leq \epsilon$ for $j = 1,\dots, d$, With Weyl's well-known equidistribution argument in mind, it is natural to take a smooth function $\chi$ which approximates the characteristic function of the cube $[-\epsilon, \epsilon]^d \in (\R/\Z)^d$ and then evaluate
\[ \sum_{n \leq N} \chi(n^2 \alpha_1,\dots, n^2 \alpha_d)\] by expanding $\chi$ as a Fourier series on $\Z^d$. Using Weyl's inequality for quadratic phases, which we will discuss shortly, such a procedure provides a good (and, in particular, positive) estimate provided that there are no ``diophantine'' relations amongst the $\alpha_j$. Problems are encountered when, for example,
\[ \Vert r_1 \alpha_1 + \dots + r_d \alpha_d \Vert_{\R/\Z} \] is small for smallish integers $r_i$. However it turns out that if there \emph{is} such a relation then it may be used to essentially reduce the dimension of the problem by one, so that one may proceed inductively.
In order to make the induction efficient one cannot work simply with cubes $[-\epsilon,\epsilon]^d$. Instead one must work with a larger class of domains, such as arbitrary symmetric convex bodies $K$. Using some arguments in the geometry of numbers or in finite-dimensional Banach space theory one may approximate $K$ by an ellipsoid $\widetilde{K}$. Thus one is interested in whether there is $n \leq N$ such that $(n^2\alpha_1,\dots,n^2 \alpha_d) \in \widetilde{K} + \Z^d$. By a linear transformation one may map $\widetilde{K}$ to the unit ball $B(0,1)$, and the problem then becomes one of determining whether $(n^2 \alpha'_1,\dots, n^2 \alpha'_d) \in B(0,1) + \Lambda$, for a lattice $\Lambda \in \R^d$. Schmidt's result says that this is so if $N$ is suitably large depending on $\det(\Lambda)$ and, as we remarked, it is essentially proved by induction on the dimension of $\Lambda$.

Our approach will be more-or-less the same. However we make the observation that a rather natural smooth approximation to the characteristic function of $B(0,1) + \Lambda$ is provided by the theta function associated to $\Lambda$. This is particularly so if one wishes to do harmonic analysis, as the Poisson summation formula takes a very pleasant form.

\begin{definition}[Theta functions]
Suppose that $\Lambda$ is a lattice of full rank in $\R^d$. For any $t > 0$ and $x = (x_1,\ldots,x_d) \in \R^d$, we define the \emph{theta function}
\[ \Theta_{\Lambda}(t,x) := \sum_{m \in \Lambda} e^{-\pi t |x-m|^2}\]
where $|x| := (x_1^2 + \ldots + x_d^2)^{1/2}$ is the usual Euclidean norm. 
\end{definition}
\begin{remark}
For most of this appendix, one should think of $\Theta_{\Lambda}(t,x)$ as a blurred version of the characteristic function of the set obtained by placing a Euclidean ball of radius $\sim 1/\sqrt{t}$ about every point of $\Lambda$.
\end{remark}

From the Poisson summation formula we have the fundamental identity
\begin{equation}\label{poisson}
\sum_{m \in \Lambda} e^{-\pi t |x-m|^2} = \frac{1}{t^{d/2} \det(\Lambda)} \sum_{\xi \in \Lambda^*} e^{-\pi |\xi|^2/t} e(\xi \cdot x)
\end{equation}
where $\Lambda^* := \{ \xi \in \R^d: \xi \cdot m \in \Z \hbox{ for all } m \in \Lambda \}$ is the dual lattice of $\Lambda$.  

The determinant $\det(\Lambda)$ is, of course, an important quantity associated with the lattice $\Lambda$. In our argument, however, a somewhat different quantity will play a more prominent r\^ole.

\begin{definition}[Definition of $A_{\Lambda}$] Let $\Lambda$ be a lattice of full rank in $\R^d$. Define
\begin{equation}\label{A-form}
A_{\Lambda} := \Theta_{\Lambda^*,d}(1,0) = \sum_{\xi \in \Lambda^*} e^{-\pi |\xi|^2} = \det(\Lambda) \sum_{m \in \Lambda} e^{-\pi |m|^2}.
\end{equation}
\end{definition}
\begin{remark} The last equality follows from \eqref{poisson}. $1/A_{\Lambda}$ may be thought of as a kind of measure of how likely it is that a random point in $\R^d$ lies within $O(1)$ of $\Lambda$.
\end{remark}

Now let $\alpha = (\alpha_1,\dots,\alpha_d) \in \R^d$ and let $N > 0$. We define the quantity
\[ F_{\Lambda,\alpha}(N) := \det(\Lambda) \E_{-N \leq n \leq N} \Theta_{\Lambda}(1, n^2 \alpha)\]
From \eqref{poisson} we have
\begin{equation}\label{f-fourier} F_{\Lambda,\alpha}(N) = \sum_{\xi \in \Lambda^*} e^{-\pi |\xi|^2} \E_{-N \leq n \leq N} e(n^2 \xi \cdot \alpha).\end{equation}We will work towards a lower bound for $F_{\Lambda,\alpha}(N)$. The precise statement of this bound may be found in Proposition \ref{prop-a9} below. Once this is available, a straightforward truncation argument can be used to show that $n^2 \alpha$ is often within $O(1)$ of the lattice $\Lambda$. Rescaling suitably, one may insist that $n^2 \alpha$ is within $\epsilon$ of $\Lambda$ under appropriate conditions, and Proposition \ref{sqj} follows. We postpone the details until the end of the section, focussing for now on the much more interesting issue of a lower bound for $F_{\Lambda,\alpha}(N)$.

Later on we will need the following list of simple but slightly technical properties of $F_{\Lambda,\alpha}$. The reader may care to skip the next lemma on a first reading.

\begin{lemma}[Properties of $F_{\Lambda,\alpha}$]\label{flambda}  Let $\Lambda$ be a lattice of full rank in $\R^d$, let $\alpha \in \R^d$, and let $N > 0$.
\begin{itemize}
\item[(i)] \textup{(Contraction of $N$)} For any $c \in (\frac{10}{N},1)$, we have $F_{\Lambda,d}(\alpha,N) \gg c F_{\Lambda,d}(\alpha, cN)$.
\item[(ii)] \textup{(Dilation of $\alpha$)} For any integer $q \geq 1$, we have $F_{\Lambda,d}(\alpha,N) \gg \frac{1}{q} F_{\Lambda,d}(q^2 \alpha, N/q)$.
\item[(iii)] \textup{(Stability)} If $\tilde \alpha \in \R^d$ is such that $|\alpha - \tilde \alpha| \leq \eps N^{-2}$ for some $\eps \in (0,1)$, 
then $F_{\Lambda,d}(\alpha,N) \gg F_{(1+\eps) \cdot \Lambda,d}((1+\eps) \tilde \alpha, N)$.
\end{itemize}
\end{lemma}
\begin{proof} The bound (i) follows immediately from the definition of $F_{\Lambda,\alpha}$, the positivity of $\Theta$. The bound (ii) also follows immediately from the definition of $F_{\Lambda,\alpha}$, restricting the $n$ variable to multiples of $q$. We now turn to the stability estimate (iii).  If $|\alpha - \tilde \alpha| \leq \eps N^{-2}$ then 
\[ \big||n^2 \alpha - m| - | n^2 \widetilde \alpha - m|\big| \leq \eps\] for all $n$, $-N \leq n \leq N$, and all $m \in \Lambda$.
Write, temporarily, $X := |n^2 \alpha - m|$ and $\widetilde{X} := |n^2 \widetilde\alpha - m|$. If $X \geq 2$ then we have the inequality
\[ \pi (1 + \eps)^2 (X - \eps)^2 \geq \pi X^2,\] and so
\[ e^{-\pi X^2} \geq e^{-(1 + \eps)^2 \pi \widetilde X^2}.\]
If $X \leq 2$ then $e^{-\pi X^2} > c$, and so 
\[ e^{-\pi X^2} \geq ce^{-(1 + \eps)^2 \pi \widetilde X^2}\] in this case.
Thus in both cases we have
 \[ e^{-\pi X^2} \gg e^{-(1+\eps)^2\pi \widetilde{X}^2}.\] Substituting for $X,\widetilde{X}$, 
summing this in $m$ and averaging in $n$, the claim follows.\end{proof}

The next lemma is the key ingredient in our argument. It formalises the idea that everything is relatively straightforward unless there is a ``diophantine'' relation amongst the $\alpha_i$.

\begin{lemma}[Schmidt's alternative]\label{alternative} Suppose that $\alpha \in R^d$ and that $\Lambda \subseteq \R^d$ is a full-rank lattice. Let $N > 0$ be an integer.
One of the following two alternatives always holds:
\begin{enumerate}
\item $F_{\Lambda,\alpha}(N) \geq 1/2$;
\item There is a positive integer $q \ll dA_{\Lambda}^C$ and some primitive $\xi \in \Lambda^* \setminus \{0\}$ such that
\begin{equation}\label{xi-size}
| \xi | \ll \sqrt{d} + \sqrt{\log A_{\Lambda}}
\end{equation}
and
\begin{equation}\label{q-bound} \Vert q\xi \cdot \alpha \Vert_{\R/\Z} \ll A_{\Lambda}^C N^{-2}.
\end{equation}
\end{enumerate}
\end{lemma}
\begin{remark}
We say that $\xi \in \Lambda^*$ is \emph{primitive} if $\xi/n \notin \Lambda^*$ for any integer $n \geq 2$.
\end{remark}
\proof Suppose that (i) fails to hold. Then from \eqref{f-fourier} and the triangle inequality we have, using \eqref{poisson}, that
\begin{equation}\label{eq8.80} \sum_{\xi \in \Lambda^* \setminus\{0\}} e^{-\pi |\xi|^2} \big| \E_{-N \leq n \leq N} e(n^2 \xi \cdot \alpha) \big| > 1/2.\end{equation}
Our first task is to truncate this. To this end let $M \geq 1$ be a cutoff parameter to be chosen later. We have
\begin{align*}
\sum_{\xi \in \Lambda^*: |\xi| \geq M} e^{-\pi |\xi|^2} |\E_{-N \leq n \leq N} e( n^2 \xi \cdot \alpha )| 
&\leq \sum_{\xi \in \Lambda^*: |\xi| \geq M} e^{-\pi |\xi|^2} \\
&\leq e^{-\pi M^2/2} \sum_{\xi \in \Lambda^*} e^{-\pi |\xi|^2/2} \\
&= e^{-\pi M^2/2} 2^{d/2} \det(\Lambda) \sum_{m \in \Lambda} e^{-2 \pi |m|^2} \\
&\leq e^{-\pi M^2/2} 2^{d/2} A_{\Lambda}
\end{align*}
Choosing $M := C(\sqrt{d} + \sqrt{\log A_{\Lambda}})$ for suitable $C$ we may clearly make this less than $1/4$, and hence from \eqref{eq8.80} we have
\[ \sum_{\xi \in \Lambda^*: 0 < |\xi| < M} e^{-\pi |\xi|^2} \big| \E_{-N \leq n \leq N} e(n^2 \xi \cdot \alpha) \big| \geq 1/4.\]
From the definition of $A_{\Lambda}$ this implies that there is $\xi \in \Lambda^* \setminus 0$, $|\xi| \leq M$, such that
\begin{equation}\label{use-weyl} | \E_{-N \leq n \leq N} e(n^2 \xi \cdot \alpha) | \geq 1/4A_{\Lambda}.\end{equation}
This puts us in the situation covered by Weyl's inequality, a discussion of which may be found in \cite[Chapter 3]{montgomery-ten} or \cite[Chapter 2]{vaughan}. The following formulation of the result follows easily from the standard one as given in those two references; see also \cite[Lemma A.13]{gt-mobius}.
\begin{weyl-theorem}
Let $\theta \in R$, let $\delta \in (0,1)$ and suppose that $N > 0$ is an integer such that $|\E_{-N \leq n \leq N}e(n^2 \theta)| \geq \delta$. Then there exists a positive integer $q \ll \delta^{-C_1}$ such that $\Vert q \theta \Vert_{\R/\Z} \ll \delta^{-C_2} N^{-2}$.
\end{weyl-theorem}
\begin{remark}
For us the exact values of $C_1,C_2$ are unimportant, but it is possible to take $C_1 = 2$ and $C_2$ to be any number larger than $2$.
\end{remark}
The bounds \eqref{xi-size} and \eqref{q-bound} follow immediately from this and \eqref{use-weyl}. It remains to show that $\xi$ can be chosen to be primitive. There is certainly a natural number $n$ such that $\xi/n$ lies in $\Lambda^*$ and is primitive. Setting $\widetilde{\xi} := \xi/n$ and $\widetilde{q} := nq$ it is clear that the bounds \eqref{xi-size} and \eqref{q-bound} are preserved. We must show that $\widetilde{q} \ll dA_{\Lambda}^{C'}$ for some absolute $C'$. To do this we note from \eqref{A-form} that if $\lambda^* \in \Lambda^* \setminus \{0\}$ is arbitrary then
\begin{equation}\label{lamstar}
A_{\Lambda} \gg 1/|\lambda^*|
\end{equation}
(consider the cases $|\lambda^*| \leq 1$ and $|\lambda^*| > 1$ separately).
It follows from this, \eqref{xi-size} and a crude bound that
\[ n = \frac{|\xi|}{|\widetilde{\xi}|} \ll A_{\Lambda}(\sqrt{d} + \sqrt{\log A_{\Lambda}}) \ll dA_{\Lambda}^2.\]
The alternative lemma follows immediately.\endproof

We will shortly combine the alternative lemma with some additional arguments which allow us to make progress when case (ii) holds. We first isolate a simple but important lemma that will be needed.

\begin{lemma}[Descent]\label{ind-help}
Suppose that $\Lambda' \subseteq \R^{d-1}$ and $\Lambda \subseteq \R^d$ are full-rank lattices, and that $\Lambda' \subseteq \Lambda$, where we are regarding $\R^{d-1}$ as a subset of $\R^d$ in the usual way. Suppose that $\alpha' \in R^{d-1}$, that $\alpha \in \R^d$ and that $\alpha - \alpha' \in \Lambda$. Then
\[ F_{\Lambda,\alpha}(N) \geq \frac{\det(\Lambda)}{\det(\Lambda')} F_{\Lambda',\alpha'}(N).\]
\end{lemma}
\proof By definition we have
\[ F_{\Lambda,\alpha}(N) = \det(\Lambda) \E_{-N \leq n \leq N} \sum_{m \in \Lambda} e^{-\pi |n^2 \alpha - m|^2},\] and there is a similar expression for $F_{\Lambda',\alpha'}$. Now by translation invariance and positivity we have, for each fixed $n$,
\[ \sum_{m \in \Lambda} e^{-\pi|n^2\alpha - m|^2} = \sum_{m \in \Lambda} e^{-\pi |n^2 \alpha' - m|^2} \geq \sum_{m \in \Lambda'} e^{-\pi |n^2 \alpha' - m'|^2}.\]
The result follows upon taking expectations over $n$.\endproof

\begin{proposition}[Inductive lower bound on $F_{\Lambda,\alpha}$]\label{iteration-step}
Suppose that $\alpha \in R^d$ and that $\Lambda \subseteq \R^d$ is a full-rank lattice. Let $N > d^CA_{\Lambda}^C$ be an integer. Then either $F_{\Lambda,\alpha}(N) \geq 1/2$ or else there is an $\alpha' \in \R^{d-1}$, a full-rank lattice $\Lambda' \subseteq \R^{d-1}$ with 
\begin{equation}\label{a-bd} A_{\Lambda'} \ll (\sqrt{d} + \sqrt{\log A_{\Lambda}}) A_{\Lambda},\end{equation}
 and an $N' \gg d^{-C}A_{\Lambda}^{-C}N$ such that
\[ F_{\Lambda,\alpha}(N) \geq d^{-C}A_{\Lambda}^{-C} F_{\Lambda',\alpha'}(N').\]\end{proposition}
\proof We begin by applying the alternative lemma. We may clearly assume that we are in case (ii), that is to say there exists a primitive $\xi \in \Lambda^* \setminus 0$ and a $q \ll d A_{\Lambda}^C$ such that \eqref{xi-size} and \eqref{q-bound} are satisfied. By subjecting $\alpha$ and $\Lambda$ to a rotation, we may assume without loss of generality that $\xi = \xi_d e_d$ is a multiple of the basis vector $e_d$. Now multiplying through by $q$ we see from \eqref{q-bound} that 
\[ \Vert \xi \cdot q^2 \alpha \Vert_{\R/\Z} \ll d A_{\Lambda}^C N^{-2}.\] Recalling \eqref{lamstar}, we can find $\beta \in \R^d$ such that $\xi \cdot \beta \in \Z$ and
\[ |\beta - q^2 \alpha| \leq |\xi_d|^{-1} \Vert \xi \cdot  q^2 \alpha \Vert_{\R/\Z} \ll dA_{\Lambda}^C N^{-2}.\]
In particular we may choose 
\[ N_* \gg d^{-C} A_{\Lambda}^{-C} N\]
such that
\begin{equation}\label{beta-bd} |\beta - q^2 \alpha| \leq  N_*^2/d.\end{equation}
Now $\xi$ is primitive, and so there is $m \in \Lambda$ so that $\xi \cdot \beta = \xi \cdot m$. Since $\xi = \xi_d e_d$ this means that we may write $\beta = \beta' + m$ where $\beta' \in \R^{d-1}$.

Now by Lemma \ref{flambda} (i) we have
\[ F_{\Lambda,\alpha}(N) \gg d^{-C} A_{\Lambda}^{-C} F_{\Lambda,\alpha}(N_*).\] Note that our lower bound on $N$ ensures that the value of $c$ in that lemma can be taken to be at least $10/N$, as required.
By Lemma \ref{flambda} (ii) and the fact that $q \ll d A_{\Lambda}^C$ this implies that
\[ F_{\Lambda,\alpha}(N) \gg d^{-C} A_{\Lambda}^{-C} F_{\Lambda,q^2 \alpha}(N_*/q).\]
Lemma \ref{flambda} (iii) and \eqref{beta-bd} allow us to assert that
\[ F_{\Lambda,\alpha}(N) \gg d^{-C} A_{\Lambda}^{-C} F_{(1 + 1/d) \Lambda, (1 + 1/d)\beta}(N_*/q).\]
From Lemma \ref{ind-help} we obtain
\begin{equation}\label{use-in-a-sec} F_{\Lambda,\alpha}(N) \gg d^{-C} A_{\Lambda}^{-C} \frac{\det(\Lambda)}{\det(\Lambda \cap \R^{d-1})} F_{\Lambda', \alpha'}(N'),\end{equation}
where $\alpha' := (1 + 1/d) \beta'$, $\Lambda' := (1 + 1/d)\Lambda \cap \R^{d-1}$ and $N' := N_*/q$. The claimed bound on $N'$ follows immediately from the lower bound on $N_*$ and the upper bound $q \ll dA_{\Lambda}^C$.

Now since $\xi$ is primitive and parallel to $e_d$ we have $\det(\Lambda^*) = |\xi| \det((\Lambda \cap \R^{d-1})^*)$. Since $\det(\Pi)\det(\Pi^*) = 1$ for any lattice $\Pi$, the ratio of determinants in \eqref{use-in-a-sec} is $|\xi|^{-1}$.In view of the upper bound \eqref{xi-size}, this may be absorbed into the $d^{-C} A_{\Lambda}^{-C}$ factor, and we therefore obtain the claimed lower bound on $F_{\Lambda,\alpha}(N)$. 

It remains to place an upper bound on $A_{\Lambda'}$. Note first that by positivity we have
\[  \frac{A_{\Lambda \cap \R^{d-1}}}{\det(\Lambda \cap \R^{d-1})} \leq  \frac{A_{\Lambda}}{\det(\Lambda)}, \] and so from the previous discussion and \eqref{xi-size} we have
\begin{equation}\label{alm} A_{\Lambda \cap \R^{d-1}} \leq |\xi|A_{\Lambda} \ll (\sqrt{d} + \sqrt{\log A_{\Lambda}}) A_{\Lambda}.\end{equation}
Secondly for any lattice $\Pi$ and any $\delta > 0$ we clearly have
\[ \sum_{m \in (1 + \delta)\Pi} e^{-\pi |m|^2} \leq \sum_{m \in \Pi} e^{-\pi |m|^2},\] 
and so
\[ A_{\Lambda'} \leq (1 + 1/d)^d A_{\Lambda \cap \R^{d-1}} \ll A_{\Lambda \cap \R^{d-1}}.\]
Combining this with \eqref{alm}, we obtain the required upper bound on $A_{\Lambda'}$.\endproof

Iterating this proposition leads in a straightforward manner to the claimed lower bound on $F_{\Lambda,\alpha}(N)$.

\begin{proposition}[Lower bound for $F_{\Lambda,\alpha}$]\label{prop-a9}
Let $\alpha \in \R^d$, suppose that $\Lambda \subseteq \R^d$ is a lattice of full rank with $\det(\Lambda) \geq 1$, and let $N > 0$ be an integer. Then we have the lower bound
\[ F_{\Lambda,\alpha}(N) \gg d^{-Cd^2}A_{\Lambda}^{-Cd}.\]
\end{proposition}
\proof If $N < d^{C_0d^2}A_{\Lambda}^{C_0d}$ then the result is immediate from the trivial lower bound
\[ F_{\Lambda,\alpha}(N) \geq \det(\Lambda)/(2N + 1).\]
Suppose then that $N \geq d^{C_0d^2}A_{\Lambda}^{C_0d}$ for some suitably large $C_0$.
Set $\alpha_0 := \alpha$, $\Lambda_0 := \Lambda$ and $N_0 := N$. Apply Proposition \ref{iteration-step} repeatedly, obtaining vectors $\alpha_j \in \R^{d-j}$, lattices $\Lambda_j \subseteq \R^{d-j}$ and integers $N_j$ for $j = 0,1,\dots$. We will show in a short while that $N_j > d^C A_{\Lambda_j}^C$ throughout this iteration, and so it is indeed valid to continue applying Proposition \ref{iteration-step}.
If, at some point, we pass through case (i) of the alternative lemma (which leads to the lower bound $F_{\Lambda,\alpha}(N) \geq 1/2$) then we stop the iteration. The worst bounds result from when this is not the case, and the iteration proceeds all the way to $d = 0$. Note that we have $F_{\Lambda,\alpha}(N) = 1$ when $d = 0$. The growth of $A_{\Lambda_j}$ during the iteration is controlled by \eqref{a-bd}. Noting that $A_{\Lambda} \geq \det(\Lambda) \geq 1$, we may employ the crude inequality
\[ \sqrt{d} + \sqrt{\log X} \ll d X^{1/d}\]
for $X \geq 1$. Using this it is easy to see\footnote{We thank Zach Hunter for drawing our attention to the fact that the published version of the paper, which omitted the dependence on $d$ here, was incorrect.} from \eqref{a-bd} that
\[ A_{\Lambda_j} \ll d^{Cj} A_{\Lambda_0}^{C}\] for the duration of the iteration.
Since
\[ N_{j+1} \geq d^{-C} A_{\Lambda_j}^{-C} N_j\] for all $j$,
 this confirms that $N_j > d^C A_{\Lambda}^C$ throughout provided that $C_0$ is chosen large enough. Since
\[ F_{\Lambda_{j},\alpha_{j}}(N_j) \gg d^{-C} A_{\Lambda_j}^{-C} F_{\Lambda_{j+1},\alpha_{j+1}}(N_{j+1}),\]
it also provides the desired lower bound on $F_{\Lambda,\alpha}(N)$.\endproof

It remains to deduce Proposition \ref{sqj}. This is achieved by a truncation argument.

\emph{Proof of Proposition \ref{sqj}.} Let $R$ be a quantity to be chosen later. We will need $R > C_0d$ for some large absolute constant $C_0$. Apply Proposition \ref{prop-a9} with $\alpha := (\alpha_1,\dots,\alpha_d)$ and $\Lambda := R\Z^d$. We have
\[ A_{\Lambda} = R^d \big( \sum_{m \in R\Z} e^{-\pi m^2} \big)^d \leq (CR)^d,\] and so (since $R \geq Cd$) that proposition implies that
\[ F_{\Lambda,\alpha}(N) \gg R^{-Cd^2}.\]
Since $\det(\Lambda) = R^d$, it follows from the definition of $F_{\Lambda,\alpha}$ that
\[ \E_{-N \leq n \leq N} \sum_{m \in R\Z^d} e^{-\pi |n^2 \alpha - m|^2} \gg R^{-Cd^2}.\]
The contribution of the $n = 0$ term is $\ll (CR)^d/N$, which is negligible if $N \geq CR^{Cd^2}$ for suitably large $C$. In this case we conlcude that there is $n \in \{1,\dots,N\}$ such that 
\begin{equation}\label{7.10} \sum_{m \in R\Z^d} e^{-\pi |n^2 \alpha - m|^2} \gg R^{-Cd^2}.\end{equation}
Fix this $n$. If we had $|n^2 \alpha - m| > \sqrt{R}$ for all $m \in R\Z^d$ then we would have
\[ e^{-\pi |n^2 \alpha - m|^2} \leq e^{-\pi R^2/2} e^{-\pi |n^2 \alpha - m|^2/2}\] for all $m \in R\Z^d$. Summing in $m$ and using \eqref{poisson} and \eqref{A-form}, we conclude that
\[ \sum_{m \in R\Z^d} e^{-\pi |n^2 \alpha - m|^2} \leq e^{-\pi R^2/2} \frac{2^{d/2}}{\det(\Lambda)} \sum_{\xi \in \Lambda^*} e^{-2\pi |\xi|^2} e(\xi \cdot n^2 \alpha) \leq e^{-\pi R^2/2} 2^{d/2} \frac{A_{\Lambda}}{\det(\Lambda)},\]
which is $\ll e^{-\pi R^2/2} (CR)^d$. Recall that $R \geq C_0 d$; if $C_0$ is chosen large enough then this will contradict \eqref{7.10}. We are thus forced to conclude that there \emph{is} some $m \in R\Z^d$ such that $|n^2 \alpha - m| \leq \sqrt{R}$, and this clearly implies that $\Vert n \alpha_j \Vert_{\R/\Z} \leq 1/\sqrt{R}$ for $j = 1,\dots,d$. 

We have shown that if $N \geq CR^{Cd^2}$ and $R \geq Cd$ then there is some $n$, $1 \leq n \leq N$, such that $\Vert n^2 \alpha_j \Vert_{\R/\Z} \leq 1/\sqrt{R}$ for $j = 1,\dots, d$. If $N \geq C' d^{C'd^2}$ for some suitably large $C'$ then the proposition follows by choosing $R = d^{-1} N^{c/d^2}$ for some small absolute constant $c > 0$; if instead $N < C' d^{C'd^2}$ then the proposition is trivial. \endproof

\end{document}